\documentclass[12pt]{amsart}
\begin{document}
\numberwithin{equation}{section}

\def\1#1{\overline{#1}}
\def\2#1{\widetilde{#1}}
\def\3#1{\widehat{#1}}
\def\4#1{\mathbb{#1}}
\def\5#1{\frak{#1}}
\def\6#1{{\mathcal{#1}}}

\def\C{{\4C}}
\def\R{{\4R}}
\def\N{{\4N}}
\def\Z{{\4Z}}

\title[New normal forms for Levi-nondegenerate hypersurfaces]{New normal forms for Levi-nondegenerate hypersurfaces}
\author[D. Zaitsev]{Dmitri Zaitsev}
\dedicatory {Dedicated to Linda Preiss Rothschild  on the occasion of her birthday}
\address{D. Zaitsev: School of Mathematics, Trinity College Dublin, Dublin 2, Ireland}
\email{zaitsev@maths.tcd.ie}

\maketitle

\def\Label#1{\label{#1}}


\def\cn{{\C^n}}
\def\cnn{{\C^{n'}}}
\def\ocn{\2{\C^n}}
\def\ocnn{\2{\C^{n'}}}


\def\dist{{\rm dist}}
\def\const{{\rm const}}
\def\rk{{\rm rank\,}}
\def\id{{\sf id}}
\def\aut{{\sf aut}}
\def\Aut{{\sf Aut}}
\def\CR{{\rm CR}}
\def\GL{{\sf GL}}
\def\Re{{\sf Re}\,}
\def\Im{{\sf Im}\,}
\def\span{\text{\rm span}}
\def\tr{{\sf tr}}

\def\codim{{\rm codim}}
\def\crd{\dim_{{\rm CR}}}
\def\crc{{\rm codim_{CR}}}

\def\phi{\varphi}
\def\eps{\varepsilon}
\def\d{\partial}
\def\a{\alpha}
\def\b{\beta}
\def\g{\gamma}
\def\G{\Gamma}
\def\D{\Delta}
\def\Om{\Omega}
\def\k{\kappa}
\def\l{\lambda}
\def\L{\Lambda}
\def\z{{\bar z}}
\def\w{{\bar w}}
\def\Z{{\1Z}}
\def\t{\tau}
\def\th{\theta}

\emergencystretch15pt
\frenchspacing

\newtheorem{Thm}{Theorem}[section]
\newtheorem{Cor}[Thm]{Corollary}
\newtheorem{Pro}[Thm]{Proposition}
\newtheorem{Lem}[Thm]{Lemma}

\theoremstyle{definition}\newtheorem{Def}[Thm]{Definition}

\theoremstyle{remark}
\newtheorem{Rem}[Thm]{Remark}
\newtheorem{Exa}[Thm]{Example}
\newtheorem{Exs}[Thm]{Examples}

\def\bl{\begin{Lem}}
\def\el{\end{Lem}}
\def\bp{\begin{Pro}}
\def\ep{\end{Pro}}
\def\bt{\begin{Thm}}
\def\et{\end{Thm}}
\def\bc{\begin{Cor}}
\def\ec{\end{Cor}}
\def\bd{\begin{Def}}
\def\ed{\end{Def}}
\def\br{\begin{Rem}}
\def\er{\end{Rem}}
\def\be{\begin{Exa}}
\def\ee{\end{Exa}}
\def\bpf{\begin{proof}}
\def\epf{\end{proof}}
\def\ben{\begin{enumerate}}
\def\een{\end{enumerate}}
\def\beq{\begin{equation}}
\def\eeq{\end{equation}}

\section{Introduction}
In this paper we construct a large class of new normal forms
for {\em Levi-nondegenerate real hypersurfaces} in complex spaces.
We adopt a general approach illustrating why these normal forms are natural
and which role is played by the celebrated Chern-Moser normal form \cite{CM}.
The latter appears in our class as the one with
the "maximum normalization" in the lowest degree.
However, there are other natural normal forms, 
even with normalization conditions for the terms of the same degree.
Some of these forms do not involve the cube of the trace operator
and, in that sense, are simplier than the one by Chern-Moser.
We have attempted to give a complete and self-contained exposition
(including proofs of well-known results about trace decompositions)
that should be accessible to graduate students.

All normal forms here are formal, i.e.\ at the level of formal power series.
This is sufficient for most purposes such as constructing invariants
or solving the local equivalence problem for real-analytic hypersurfaces. 
In fact, a formal equivalence map between Levi-nondegenerate real-analytic
hypersurfaces is automatically convergent and is therefore a local biholomorphic map.
This is a special case of an important result of Baouendi-Ebenfelt-Rothschild
\cite{BERjams} and can also be obtained from the Chern-Moser theory \cite{CM}.
The reader is referred for more details to an excellent survey \cite{BERbull}.
We also refer to normal forms for {\em Levi-degenerate hypersurfaces}
\cite{Eindiana, Ejdg, Ko}, for real CR submanifolds of {\em higher codimension} \cite{ES,SS},
for submanifolds with {\em CR singularities} \cite{MW, HY},
and for {\em non-integrable Levi-nondegenerate hypersurface type CR structures} \cite{Z}.

Throughout the paper we consider a real-analytic hypersurface $M$ in $\C^{n+1}$
passing through $0$ and locally given by an equation
\begin{equation}\Label{graph}
\Im w=\phi(z,\bar z,\Re w), \quad (z,w)\in\C^n\times\C,
\end{equation}
where we think of $\phi$ as a power series in the components of $z,\bar z\in\C^n$ and 
$u=\Re w\in \R$. If the given hypersurface is merely smooth, one can still consider
the Taylor series of its defining equation, which is a formal power series.
This motivates the notion of  a {\em formal hypersurface}, i.e.\ the one
given (in suitable coordinates) by \eqref{graph} with $\phi$ being a formal power series.

In order to describe the normal forms, consider the expansion
\begin{equation}\Label{expand-phi}
\phi(z,\bar z,u)= \sum_{k,l,m} \phi_{kml}(z,\bar z)u^l,
\end{equation}
where each $\phi_{kml}(z,\bar z)$ is a bihomogeneous polynomial in $(z,\bar z)$
of bidegree $(k,m)$, i.e.\ $\phi_{kml}(tz,s\bar z) = t^k s^m \phi_{kml}(z,\bar z)$ for $t,s\in\R$.
The property that $\phi$ is real is equivalent to the reality conditions
\begin{equation}\Label{reality}
\1{\phi_{kml}(z,\bar z)} = \phi_{mkl}(z,\bar z).
\end{equation}
The property that $M$ passes through $0$ corresponds to $\phi$ having {\em no constant terms}.
Furthermore, after a complex-linear transformation of $\C^{n+1}$,
one may assume that $\phi$ also has {\em no linear terms},
which will be our assumption from now on.
The {\em Levi form} now corresponds to $\phi_{11}(z,\bar z)$,
the only lowest order term that cannot be eliminated after a biholomorphic change of coordinates. Following \cite{CM}, we write $\phi_{11}(z,\bar z)=\langle z,z\rangle$,
which is a hermitian form in view of \eqref{reality}.
If the Levi form is nondegenerate,
after a further complex-linear change of coordinates, we may assume that
\begin{equation}\Label{levi-exp}
\langle z, z\rangle = \sum_{j=1}^n \eps_j z_j \bar z_j, \quad \eps_j=\pm1.
\end{equation}
An important role in the normal forms is played by the 
{\em trace operator} associated with $\langle z, z\rangle$, which is
the second order differential operator given by
\begin{equation}\Label{trace}
\tr:= \sum_{j=1}^n \eps_j \frac{\d^2}{\d z_j \d \bar z_j}.
\end{equation}
In particular, for $n=1$ we have 
\begin{equation}\Label{n2}
\phi_{kml}(z,\bar z)=c_{kml}z^k\bar z^m, \quad 
\tr\phi_{kml}=\tr (c_{kml}z^k\bar z^m) =
\begin{cases}
km c_{kml}z^{k-1}\bar z^{m-1}, & \min(k,m)\ge 1,\\
0, & \min(k,m)=0.
\end{cases}
\end{equation}

We now consider the following normalization conditions.
\begin{enumerate}
\item For every $k\ge2$ and $l\ge1$, choose disjoint integers $0\le m,m'\le l$ with $m\ge1$, $m'\ne m$
and consider the conditions
\begin{equation*}
\tr^{m-1}\phi_{k+m,m,l-m}=0,\quad \tr^{m'}\phi_{k+m',m',l-m'}=0.
\end{equation*}
\item For every $k\ge 2$, consider the condition
\begin{equation*}
\phi_{k00}=0.
\end{equation*}
\item For every $l\ge2$, choose pairwise disjoint integers $0\le m,m',m''\le l$ with $m\ge 1$,
of which the nonzero ones are not of the same parity (i.e.\ neither all are even nor all are odd)
and consider the conditions
\begin{equation*}
\tr^{m-1}\phi_{m+1,m,l-m}=0, \quad \tr^{m'}\phi_{m'+1,m',l-m'}=0, \quad \tr^{m''}\phi_{m''+1,m'',l-m''}=0.
\end{equation*}
\item For every $l\ge 3$, choose disjoint even integers $0\le m,m'\le l$ with $m\ge1$
and disjoint odd integers $0\le \2m,\2m'\le l$ and consider the conditions
\begin{equation*}
\tr^{m-1}\phi_{m,m,l-m}=0, \quad  \tr^{m'}\phi_{m',m',l-m'}=0,\quad
\tr^{\2m-1}\phi_{\2m,\2m,l-\2m}=0, \quad  \tr^{\2m'}\phi_{\2m',\2m',l-\2m'}=0.
\end{equation*}
\item Consider the conditions
\begin{equation*}
\phi_{101}=0, \quad \phi_{210}=0, 
\quad \phi_{002}=0, \quad \phi_{111}=0, \quad \tr\phi_{220}=0.
\end{equation*}
\end{enumerate}

In view of \eqref{n2}, in case $n=1$ all traces can be omitted.
The multi-indices $(k,m,l)$ of $\phi_{kml}$ involved in different conditions (1)--(5)
are all disjoint. They are located on parallell lines in the direction of  the vector $(1,1,-1)$.
In fact, (1) corresponds to the lines through $(k+l,l,0)$ with $k\ge2$, $l\ge1$,
whereas (2) corresponds to the same lines with $l=0$ containing only one
triple $(k,0,0)$ with nonnegative components.
Condition (3) corresponds to the lines through $(k+l,l,0)$ for $k=1$, $l\ge 2$,
whereas (4) corresponds to the same lines for $k=0$ and $l\ge3$.
Finally (5) involves all $5$ coefficients that correspond
to the lines through $(2,1,0)$ and $(2,2,0)$.

The following is the main result concerning the above normal forms.
\bt\Label{main}
Every (formal) real hypersurface $M$ through $0$ admits a local (formal) biholomorphic transformation $h$ preserving $0$ into each of the normal forms given by (1)--(5).
If $M$ is of the form \eqref{graph} with $\phi$ having no constant and linear terms
and satisfying \eqref{levi-exp}, the corresponding transformation $h=(f,g)$
is unique provided it is normalized as follows:
\begin{equation}\Label{fg-norm}
f_{z}=\id, \quad f_w=0, \quad g_{z}=0, \quad g_{w}=1, \quad \Re g_{z^2}=0,
\end{equation}
where the subscripts denote the derivatives taken at the origin.
\et

The Chern-Moser normal form corresponds to a choice
of the coefficients of the lowest degree that may appear in (1) -- (5).
In fact, we choose $m=1$, $m'=0$ in (1), $m=1$, $m'=0$, $m''=2$ in (2),
and $m=2$, $m'=0$, $\2m=1$, $\2m'=3$ in (4) and obtain the familiar Chern-Moser normal form
\begin{equation}\Label{normal-cm}
\phi_{k0l}=0, \quad \phi_{k1l}=0, \quad \tr \phi_{22l}=0,  \quad \tr^2 \phi_{32l}=0,  
\quad \tr^3 \phi_{33l}=0,
\end{equation}
where $(k,l)\ne (1,0)$ in the second equation.

However, we also obtain other normal forms involving the same 
coefficients $\phi_{kml}$. Namely, we can exchange
$\2m$ and $\2m'$ in (4), i.e.\ choose $\2m=3$, $\2m'=1$, which leads to the normalization
\begin{equation}\Label{normal-1}
\phi_{k0l}=0, \quad \phi_{k1l}=0, \quad \tr\phi_{11l}=0, \quad \phi_{111}=0, \quad \tr \phi_{22l}=0,  \quad \tr^2 \phi_{32l}=0,  \quad \tr^2 \phi_{33l}=0,
\end{equation}
where $k\ge2$ in the second and $l\ge2$ in the third equations.
This normal form is simplier than \eqref{normal-cm}
in the sense that it only involves the trace operator and its square
rather than its cube. 
A comparison of \eqref{normal-cm} and \eqref{normal-1}
shows that the first set of conditions has more equations for $\phi_{11l}$,
whereas the second set has more equations for $\phi_{33l}$.
Thus we can say that the Chern-Moser normal form 
has the "maximum normalization" in the lowest degree.

Alternatively, we can exchange $m$ and $m''$ in (2),
i.e.\ choose $m=2$, $m'=0$, $m''=1$, 
which leads to the normalization
\begin{equation}\Label{normal-2}
\phi_{k0l}=0, \quad \phi_{k1l}=0, \quad \tr\phi_{21l}=0,
\quad \phi_{210}=0,
\quad \tr \phi_{22l}=0,  \quad \tr \phi_{32l}=0,  
\quad \tr^3 \phi_{33l}=0,
\end{equation}
where $k=1$, $l\ge 1$ or $k\ge3$ in the second and $l\ge1$ in the third equations.
Again, comparing with \eqref{normal-cm},
we see that the latter has more equations for $\phi_{21l}$
and less for $\phi_{32l}$, i.e.\ again the Chern-Moser normal form 
has the "maximum normalization" in the lowest degree.

Finally, we can combine both changes leading to \eqref{normal-1} and \eqref{normal-2},
i.e.\ choose $m=2$, $m'=0$, $m''=1$ in (2) and $\2m=3$, $\2m'=1$ in (4)
and obtain yet another normal form involving the same terms $\phi_{kml}$.
We leave it to the reader to write the explicit equations.

Thus we have 4 normal forms involving the same terms as the one by Chern-Moser.
We now describe a completely different normal form, which also has certain extremality 
property. Roughly speaking, in those 4 normal forms,
we have chosen the multi-indices with the smallest first two components on each line of indices in (1), (3), (4). We now choose the multi-indices with the smallest last component, i.e.\ those at the {\em other end} of each line. For instance, we can 
choose $m=l$, $m'=l-1$ in (1), $m=l$, $m'=l-1$, $m''=l-2$ in (3)
and $m=l$, $m'=l-2$, $\2m=l-1$, $\2m'=l-3$ 
or $\2m=l$, $\2m'=l-2$, $m=l-1$, $m'=l-3$ in (4), depending on the parity of $l$.
We obtain
\begin{equation}\Label{nu}
\tr^{l-1}\phi_{kl0}=0, \quad \phi_{k00}=0, \quad \tr^l\phi_{kl1}=0, 
\quad \tr^{l-1} \phi_{ll1}=0,  \quad \tr^l \phi_{l+1,l,2}=0,  
\quad \tr^l\phi_{ll2}=0,
\quad \tr^l \phi_{ll3}=0,
\end{equation}
where $k\ge l\ge1$, $(k,l)\ne(1,1)$ in the first, $k\ge l+1$ in the third, and $l\ge1$ in the forth equations. This normal form is distinguished by the property 
that it has the "maximum normalization" of the terms $\phi_{kml}$ with the lowest index $l$.
In fact, \eqref{nu} only involves $\phi_{kml}$ with $l=0,1,2,3$.
In particular, most harmonic terms $\phi_{k0l}$ are not eliminated in contrast to 
the normal forms \eqref{normal-cm}--\eqref{normal-2}.
In case $n=2$ all traces can be removed, in particular,
all terms $\phi_{kml}$ with $l=0$, except the Levi form, are eliminated, i.e.\
$\phi(z,\bar z,u)=\langle z,z\rangle + O(|u|)$.

Finally, we consider another interesting normal form that
in some sense mixes the one by Chern-Moser with the one in \eqref{nu}
(or, more precisely, the ones in \eqref{normal-1} and \eqref{nu}).
Here the multi-indices on both ends of the lines are involved.
This normal form corresponds to the choice of 
$m=l$, $m'=0$ in (1), $m=l$, $m'=0$, $m''=l-1$ in (3)
and $m=l$, $m'=1$, $\2m=l-1$, $\2m'=0$ if $l$ is even
and $\2m=l$, $\2m'=1$, $m=l-1$, $m'=0$ if $l$ is odd.
This leads to the normalization conditions
\begin{equation}\Label{n}
\phi_{k0l}=0, \quad \tr^{l-1}\phi_{kl0}=0, \quad \tr^l\phi_{l+1,l,1}=0, 
\quad \tr^{l-1} \phi_{ll1}=0,  \quad \tr \phi_{11l}=0,  
\end{equation}
where $k\ge l\ge 1$, $(k,l)\ne(1,1)$ in the second, $l\ge1$ in the forth,
and $l\ge 2$ in the fifth equation.
A remarkable feature of this normal form that distinguishes it from
the previous ones, including the one by Chern-Moser,
is that it only involves $\phi_{kml}$ with $\min(k,m,l)\le 1$.

The rest of the paper is now devoted to the proof of Theorem~\ref{main},
where we also explain how the normalization conditions (1)--(5) arise
and why they are natural.

\section{Transformation rule and its expansion}\Label{harmonic}
As before, we shall consider a real-analytic hypersurface $M$ in $\C^{n+1}$ through $0$,
locally given by \eqref{graph} with $\phi$ having no constant and linear terms.
To the hypersurface $M$ (or more precisely to the germ $(M,0)$)
we apply a local biholomorphic transformation
$(z,w)\mapsto (f(z,w),g(z,w))$ preserving $0$ and
transforming it into another germ $(M',0)$ of a real-analytic hypersurface in $\C^{n+1}$,
still given by an equation
\begin{equation}
\Im w'=\phi'(z',\bar z',\Re w'),
\end{equation}
where $\phi'$ has no constant and linear terms.
We consider (multi)homogeneous power series expansions
\begin{equation}\Label{expansions}
\begin{split}
&f(z,w)=\sum f_{kl}(z)w^l, \quad
g(z,w)=\sum g_{kl}(z)w^l, \\
&\phi(z,\bar z,u)= \sum \phi_{kml}(z,\bar z)u^l, \quad
\phi'(z',\bar z',u')= \sum \phi'_{kml}(z',\bar z')u'^l,
\end{split}
\end{equation}
where $f_{kl}(z)$ and $g_{kl}(z)$
are homogeneous polynomials in $z\in\C^n$ of degree $k$
and $\phi_{kml}(z,\bar z)$ and $\phi'_{kml}(z',\bar z')$ are bihomogeneous polynomials 
in $(z,\bar z)\in \C^n\times\C^n$ and $(z',\bar z')\in \C^n\times\C^n$ respectively
of bidegree $(k,l)$. Furthermore, $f_{kl}$ and $g_{kl}$ are abitrary
whereas $\phi_{kml}(z,\bar z)$ and $\phi'_{kml}(z',\bar z')$ satisfies the reality condition 
\eqref{reality}
which are equivalent to $\phi$ and $\phi'$ being real-valued.

The fact that the map $(z,w)\mapsto(f(z,w),g(z,w))$ 
transforms $(M,0)$ into $(M',0)$ can be expressed by the equation
\begin{equation}\Label{basic-equation}
\Im g(z,u+i\phi(z,\bar z,u)) = \phi'\big(f(z,u+i\phi(z,\bar z,u)), \1{f(z,u+i\phi(z,\bar z,u))},
\Re g(z,u+i\phi(z,\bar z,u))\big).
\end{equation}
We use \eqref{expansions} to expand both sides of \eqref{basic-equation}:

\begin{equation}\Label{g-exp}
\Im g(z,u+i\phi(z,\bar z,u)) = \Im\Big(
\sum g_{kl}(z)\big(u+i\sum \phi_{jhm}(z,\bar z)u^m\big)^l
\Big),
\end{equation}
\begin{multline}\Label{phi-exp}
\phi'\big(f(z,u+i\phi(z,\bar z,u)), \1{f(z,u+i\phi(z,\bar z,u))},
\Re g(z,u+i\phi(z,\bar z,u))\big) =\\
\sum \phi'_{klm}\Big(
\sum f_{ab}(z)\big(u+i\sum \phi_{jhr}(z,\bar z)u^r\big)^b,
\1{\sum f_{cd}(z)\big(u+i\sum \phi_{jhr}(z,\bar z)u^r\big)^d}
\Big) \times \\
\Big(\Re\big(\sum g_{st}(z)\big(u+i\sum \phi_{jhr}(z,\bar z)u^r\big)^t\big)\Big)^m.
\end{multline}

Since both $\phi$ and $\phi'$ have vanishing linear terms,
collecting the linear terms in \eqref{g-exp} and \eqref{phi-exp}
and substituting into \eqref{basic-equation} we obtain
\begin{equation}\Label{g-vanish}
g_{10}=0, \quad \Im g_{01} =0.
\end{equation}
Conditions \eqref{g-vanish} express the fact
that the map $(f,g)$ sends the normalized tangent space $T_0M=\C^n_z\times\R_u$
into the normalized tangent space $T_0M'=\C^n_{z'}\times\R_{u'}$ (where $u'=\Re w'$).
The first condition in \eqref{g-vanish}
implies that the Jacobian matrix of $(f,g)$ at $0$ is block-triangular
and hence its invertibility is equivalent to the invertibility of both diagonal blocks 
$f_{10}=f_z$ and $g_{01}=g_w$.

\section{Partial normalization in general case}
Our first goal is to obtain a general normalization procedure
that works for all series $\phi$ without any nondegeneracy assumption. 
The key starting point consists of identifying 
non-vanishing factors in \eqref{g-exp} and \eqref{phi-exp}.
These are $f_{10}$ and $g_{01}$. All other factors may vanish.
Then we look for terms in the expansions of \eqref{g-exp} and \eqref{phi-exp} 
involving at most one factor that may vanish. 
These are
\begin{equation}\Label{good-terms}
\Im g_{kl}(z)u^l, \quad g_{01}\phi_{kml}(z,\bar z)u^l, \quad \phi'_{kml}(f_{10}(z),\1{f_{10}(z)})(g_{01} u)^l,
\end{equation}
where we used the reality of  $g_{01}$ and \eqref{reality} and have dropped the argument $z$  for $g_{01}$ since the latter is a constant. 
The first term in \eqref{good-terms} has only $g_{kl}$ which may vanish, the second $\phi_{klm}$ and the third 
$\phi'_{klm}$. Other summands have more than one entry (term) that may vanish.
The terms \eqref{good-terms} play a crucial role in the normalization and are called here
the ``good" terms. The other terms are called the ``bad'' terms.``Good" terms can be used to obtain a partial normalization of $M$ as follows.

Consider in \eqref{basic-equation} the terms of multi-degree $(k,0,l)$ in $(z,\bar z,u)$,
where the first term from \eqref{good-terms} appears. We obtain
\begin{equation}\Label{dots-eq}
\begin{split}
\frac{1}{2i}g_{kl}(z)u^l + g_{01}\phi_{k0l}(z,\bar z)u^l &= \phi'_{k0l}(f_{10}(z),\1{f_{10}(z)})(g_{01}u)^l + \ldots, \quad k>0,\\
\Im g_{0l}u^l + g_{01}\phi_{00l}u^l &= \phi'_{00l}(g_{01}u)^l + \ldots,
\end{split}
\end{equation}
where the dots stand for the ``bad'' terms.
If no dots were present, one can suitably choose $g_{kl}$
in the first equation and $\Im g_{0l}$ in the second 
to obtain the normalization 
\begin{equation}\Label{part-norm}
\phi'_{k0l}(z)=0.
\end{equation}
In presence of the ``bad'' terms, an induction argument is used.
In fact, an inspection of the expansions  of \eqref{g-exp} and \eqref{phi-exp} 
shows that
 the terms in \eqref{dots-eq} included in the dots in \eqref{dots-eq}, involve other coefficients $g_{st}(z)$
only of order $s+t$ less than $k+l$. Thus the normalization \eqref{part-norm} can be obtained
by induction on the order $k+l$. Furthermore, the expansion terms $g_{kl}(z)$ for $k>0$
and $\Im g_{0l}$
are uniquely determined by this normalization.
However, the infinitely many terms $f_{kl}(z)$ and the real parts $\Re g_{0l}$ are not determined and act as free parameters. For every choice of those parameters,
one obtains a germ $(M',0)$ which is biholomorphically equivalent to $(M,0)$
and satisfies the normalization \eqref{part-norm}.
Thus this normalization is ``partial''.
Normalization \eqref{part-norm} is the well-known elimination of the so-called {\em harmonic terms} and works along the same lines also when $M$ is of higher codimension,
i.e.\ when $w$  is vector. 

\section{Levi-nondegenerate case}
The ``good" terms in \eqref{good-terms} were not enough to determine $f_{kl}(z)$
and $\Re g_{0l}$.
Thus we need more ``good" terms for a complete normal form.
These new ``good'' terms must be obtained from the expansion of \eqref{phi-exp}
since  $f_{kl}(z)$ do not appear in \eqref{g-exp}.
Thus we need a nonvanishing property for some $\phi'_{kml}$.
The well-known lowest order invariant of $(M,0)$ is the Levi form $\phi'_{110}(z,\bar z)$
which we shall write following \cite{CM} as $\langle z,z\rangle$.
In view of \eqref{reality}, $\langle z,z\rangle$ is a hermitian form.
We assume it to be complex-linear in the first and complex-antilinear in the second argument.
The basic assumption is now that the Levi form is {\em nondegenerate}.

Once the class of all $M$ is restricted to Levi-nondgenerate ones, we obtain 
further "good" terms involving the new nonvanishing factor $\langle z,z\rangle$. The ``good" terms in both \eqref{phi-exp} and \eqref{g-exp} come now
from the expansion of
\begin{equation}\Label{good-terms'}
\begin{split}
&\Im\big( g_{kl}(z)(u+i\langle z,z\rangle)^l\big), \quad 2\Re\big(\langle f_{kl}(z), f_{10}(z)\rangle (u+i\langle z,z\rangle)^l\big),\\
&g_{01}\phi_{kml}(z,\bar z)u^l, \quad \phi'_{kml}(f_{10}(z),\1{f_{10}(z)})(g_{01} u)^l.
\end{split}
\end{equation}
This time each coefficients $g_{kl}$ and $f_{kl}$ appears in some ``good" term
and thus can be potentially uniquely determined.
The linear coefficients $g_{01}$ and $f_{10}(z)$ play a special role.
They appear by themselves in the ``good" terms $\Im g_{01}$, $\Re g_{01}\langle z, z\rangle$
and $\langle f_{10}(z), f_{10}(z)\rangle$ in
the expansions of  \eqref{basic-equation} in multidegrees $(0,0,1)$ and $(1,1,0)$.
Equating terms of those multidegrees, we obtain
\begin{equation}\Label{linear-cond}
\Im g_{01} =0, \quad \Re g_{01}\langle z, z\rangle = \langle f_{10}(z), f_{10}(z)\rangle,
\end{equation}
where the first condition already appeared in \eqref{g-vanish}
and the second expresses the invariance of the Levi form (i.e.\ its transformation
rule as tensor).
The restrictions \eqref{linear-cond} describe all possible values of $g_{01}$ and $f_{10}(z)$,
which form precisely the group $G_0$ of all linear automorphism of the hyperquadric
$\Im w = \langle z,z\rangle$. In order to study the action by more general biholomorphic
maps $(z,w)\mapsto (f(z,w),g(z,w))$ satisfying \eqref{linear-cond}, it is convenient to write a general map as a composition of one from $G_0$ and one satisfying
\begin{equation}\Label{map-normal}
g_{01}=1, \quad f_{10}(z)=z,
\end{equation}
and study their actions separately.
Since the action by the linear group $G_0$ is easy,
we shall in the sequel consider maps $(f,g)$ satisfying \eqref{map-normal},
unless specified otherwise.

Since the Levi form $\langle z,z\rangle$ is of bidegree $(1,1)$ in $(z,\bar z)$,
we conclude from the binomial expansion of the powers in the first line of \eqref{good-terms'} that every coefficient $g_{kl}$ contributes to the ``good" terms in multidegrees
$(k+m,m,l-m)$ and $(m,k+m,l-m)$ in $(z,\bar z,u)$ for all possible $0\le m\le l$.
The latters are the integral points with nonnegative components of the lines passing through 
the points $(k,0,l)$ and $(0,k,l)$ in the direction $(1,1,-1)$.
Similarly, every coefficient $f_{kl}$ contributes to the multidegrees
$(k+m,m+1,l-m)$ and $(m+1,k+m,l-m)$ for all possible $0\le m\le l$,
corresponding to the lines passing through $(k,1,l)$ and $(1,k,l)$ in the same direction $(1,1,-1)$.
For the convenience, we shall allow one or both of $k,l$
being negative, in which case the corresponding terms are assumed to be zero.
We see that ``good" terms with $g_{kl}$ for $k>0$ appear in two lines,
whereas for $g_{0l}$ both lines coincide with that through $(0,0,l)$.
Similarly, each $f_{kl}$ with $k>1$ or $k=0$ appears in two lines,
whereas for $f_{1l}$ both lines coincide with that through $(1,1,l)$.
Thus the lines through $(0,0,l)$ are special as well as the lines
through $(1,0,l)$ and $(0,1,l)$ next to it. The latter lines contain ``good" terms
with $f_{0l}$, $f_{2,l-1}$ and $g_{1l}$. Each other line contains ``good" terms
with precisely one $f_{kl}$ and one $g_{kl}$.

Thus we treat those groups of lines separately.
Collecting in \eqref{basic-equation} terms of multi-degree $(k+m,m,l-m)$ in $(z,\bar z,u)$ for $k\ge 2$ we obtain 
\begin{multline}\Label{dots-eq1}
\frac{1}{2i}{l\choose m} g_{kl}(z)u^{l-m}(i\langle z,z\rangle)^m + \phi_{k+m,m,l-m}(z,\bar z)u^{l-m} \\
= {l-1\choose m-1}\langle f_{k+1,l-1}(z), z\rangle u^{l-m}(i\langle z,z\rangle)^{m-1}
+\phi'_{k+m,m,l-m}(z,\bar z)u^{l-m} + \ldots, 
\end{multline}
where as before the dots stand for all ``bad" terms.
Note that due to our convention, for $m=0$, the term with ${l-1\choose m-1}=0$
is not present.
Similarly we collect terms of multi-degree $(m+1,m,l-m)$,
this time we obtain
two different terms with $f_{ab}$:
\begin{multline}\Label{dots-eq2}
\frac{1}{2i}{l\choose m} g_{1l}(z)u^{l-m}(i\langle z,z\rangle)^m + \phi_{m+1,m,l-m}(z,\bar z)u^{l-m} \\
= {l-1\choose m-1}\langle f_{2,l-1}(z), z\rangle u^{l-m}(i\langle z,z\rangle)^{m-1}
+{l-1\choose m-1}\langle z, f_{0l}\rangle u^{l-m} (-i\langle z,z\rangle)^m\\
+\phi'_{m+1,m,l-m}(z,\bar z)u^{l-m} + \ldots, 
\end{multline}
where we have dropped the argument $z$ from $f_{0l}$ since the latter is a constant.
Finally, for the terms of multidegree $(m,m,l-m)$, we have
\begin{multline}\Label{dots-eq3}
{l\choose m}\Im \big(g_{0l}u^{l-m}(i\langle z,z\rangle)^m\big) + \phi_{m,m,l-m}(z,\bar z)u^{l-m} \\
= 2{l-1\choose m-1}\Re\big(\langle f_{1,l-1}(z), z\rangle u^{l-m}(i\langle z,z\rangle)^{m-1}\big)+\phi'_{m,m,l-m}(z,\bar z)u^{l-m} + \ldots.
\end{multline}

Since both sides of  \eqref{basic-equation} are real,
its multihomogeneous part of a multi-degree $(a,b,c)$
is conjugate to that of multi-degree $(b,a,c)$.
Hence the system of all equations in \eqref{dots-eq1} -- \eqref{dots-eq3}
is equivalent to \eqref{basic-equation},
i.e.\ to the property that the map $(z,w)\mapsto (f(z,w),g(z,w))$
sends $(M,0)$ into $(M',0)$.

\section{Weight estimates}\Label{weight}
In order to set up induction as in \S\ref{harmonic},
we have to estimate the degrees of $f_{kl}$ and $g_{kl}$ appearing
in the ``bad" terms in \eqref{dots-eq1} -- \eqref{dots-eq3}
and compare it with the degrees of the ``good" terms.
However, different ``good" terms with $g_{kl}$ (see \eqref{good-terms'}) do not have the same degree but rather have the same {\em weight} $k+2l$, where the weight of $z$ and $\bar z$ is $1$ and the weight of $u$ is $2$. Similarly, different ``good" terms with $f_{kl}$ have the same weight $k+2l+1$. Hence this weight is more suitable for the needed estimate.
 
We now inspect the weights of the ``bad" terms in the expansions of \eqref{g-exp}
and \eqref{phi-exp}. Recall that both $\phi$ and $\phi'$ have no constant or linear terms.
Hence the weight of $\phi_{jhm}(z,\bar z)u^m$ is greater than $2$ unless
$(j,h,m)\in \{ (1,1,0), (2,0,0), (0,2,0)\}$. In particular, the expansion of 
$$g_{kl}(z)(u+i\langle z, z\rangle + i\phi_{200}(z,\bar z) + i\phi_{020}(z,\bar z)  )^l $$
contains ``bad" terms of the same weight $k+2l$ as "good" terms.
The latter fact is not suitable for setting up an induction as in \S\ref{harmonic}.
This problem is solved by initial prenormalization of $(M,0)$ as follows.
According to \S\ref{harmonic}, one can always eliminate harmonic terms
from the expansion of $\phi$. For our purposes, it will suffice to eliminate
$\phi_{200}$ (and hence $\phi_{020}$ in view of \eqref{reality}).
Thus in the sequel, we shall assume that $\phi_{200},\phi_{020},\phi'_{200},\phi'_{020}$
are all zero.

With that assumption in mind, coming back to \eqref{g-exp},
we see that the weight of $\phi_{jhm}(z,\bar z)u^m$ always greater than $2$
unless $(j,h,m)=(1,1,0)$, i.e.\ $\phi_{jhm}(z,\bar z)u^m=\langle z, z\rangle$.
Since any ``bad" term in the expansion of 
$g_{kl}(z)\big(u+i\sum \phi_{jhm}(z,\bar z)u^m\big)^l$
contains at least one factor $\phi_{jhm}(z,\bar z)u^m$ with $(j,h,m)\ne(1,1,0)$,
its weight is greater than $k+2l$. Since the weights of all terms  (including ``bad" ones) 
in \eqref{dots-eq1} -- \eqref{dots-eq3} are $k+2l$, the ``bad" terms there coming from \eqref{g-exp}
can only contain $g_{ab}(z)$ or $\1{g_{ab}(z)}$ with weight $a+2b<k+2l$, which is now suitable for our induction.
Inspecting now the terms with $g_{st}(z)$ and $\1{g_{st}(z)}$ in the expansion of \eqref{phi-exp}, 
we see that their weights must be greater than $s+2t$.
Hence, the ``bad" terms in \eqref{dots-eq1} -- \eqref{dots-eq3} coming from \eqref{phi-exp}
can only contain $g_{ab}$ or its conjugate with weight $a+2b<k+2l$.

Similarly, we inspect ``bad'' terms containing $f_{ab}$.
This time we only need to look at the expansion of \eqref{phi-exp}.
A term in the expansion of 
\begin{equation}\Label{f-exp}
f_{ab}(z)\big(u+i\sum \phi_{jhm}(z,\bar z)u^m\big)^b
\end{equation}
is of weight greater than $a+2b$ unless it appears in the expansion of 
\begin{equation}\Label{f-exp'}
f_{ab}(z)\big(u+i\langle z,\bar z\rangle\big)^b.
\end{equation}
Keeping in mind that $\phi'$ has no linear terms, we conclude
that a ``bad" term in the expansion of \eqref{phi-exp} containing $f_{ab}$
is always of weight greater than $a+2b+1$.
The same holds for ``bad" terms containing $\1{f_{ab}}$.
Hence, a ``bad" term in \eqref{dots-eq1} -- \eqref{dots-eq3}
can only contain $f_{ab}$ or its conjugate with weight $a+2b<k+2l-1$.

Summarizing, we obtain that ``bad" terms in \eqref{dots-eq1} -- \eqref{dots-eq3}
can only contain $g_{ab}$ or its conjugate with weight $a+2b<k+2l$
and $f_{ab}$ or its conjugate with weight $a+2b<k+2l-1$.
On the other hand, the ``good" terms contain $g_{ab}$ (or $\Im g_{ab}$)
of weight precisely $k+2l$ and $f_{ab}$ of weight precisely $2k+l-1$.
Thus we may assume by induction on the weight
that all terms denoted by dots in \eqref{dots-eq1} -- \eqref{dots-eq3} are fixed
and proceed by normalizing the ``good" terms there.

\section{Trace decompositions}
The ``good" terms involving $g_{kl}(z)$
appear as products of the latters and a power of the Levi form $\langle z, z\rangle$.
As $g_{kl}(z)$ varies, these products 
\begin{equation}\Label{product}
g_{kl}(z)\langle z, z\rangle^s
\end{equation}
form a vector subspace of the space of all bihomogeneous polynomials in $(z,\bar z)$
of the corresponding bidigree $(k+s,s)$.
In order to normalize such a product, we need
to construct a complementary space to the space of all products \eqref{product}.
The latter is done by using the so-called trace decompositions described as follows.

Since the Levi form $\langle z, z\rangle$ is assumed to be nondegenerate,
we can choose coordinates $z=(z_1,\dots,z_n)$ such that \eqref{levi-exp} is satisfied
and consider the associated trace operator \eqref{trace}.
The rest of this section is devoted to the proof of the following
well-known trace decompositions (see \cite{F,S,ES} for more general decomposition results):

\bp\Label{main-decomp}
For every polynomial $P(z,\bar z)$, there exist unique polynomials
$Q(z,\bar z)$ and $R(z,\bar z)$ such that
\begin{equation}\Label{decomp}
P(z,\bar z)=Q(z,\bar z)\langle z, z\rangle^s + R(z,\bar z), \quad \tr^s R=0.
\end{equation}
\ep
Taking bihomogeneous components of all terms in \eqref{decomp}
and using the uniqueness we obtain:
\bl
If $P$ in Proposition~\ref{main-decomp} is bihomogeneous in $(z,\bar z)$, so are $Q$ and $R$.
\el
Since the trace operator is real (maps real functions into real ones),
we can take real parts of both sides in \eqref{decomp} and use the uniqueness to obtain:
\bl
If $P$ in Proposition~\ref{main-decomp} is real, so are $Q$ and $R$.
\el

We begin the proof of Proposition~\ref{main-decomp} with the following elementary lemma.

\bl
Let $P(z,\bar z)$ be a bihomogeneous polynomial of bidegree $(p,q)$. Then
\begin{equation}\Label{hom}
\sum_j P_{z_j}(z,\bar z) z_j  = p P(z,\bar z), 
\quad 
\sum_j P_{\bar z_j}(z,\bar z) \bar z_j  = q P(z,\bar z), 
\end{equation}
\el
\bpf
By the assumption, $P(sz,t\bar z)= s^p t^q P(z,\bar z)$.
Differentiating in $s$ for $s=t=1$ we obtain the first identity in \eqref{hom}.
Similarly, differentiating in $t$ for $s=t=1$ we obtain the second identity.
\epf

The following is the key lemma in the proof of Proposition~\ref{main-decomp}.

\bl
Let $P(z,\bar z)$ be a bihomogeneous polynomial of bidegree $(p,q)$. Then
\begin{equation}\Label{tr-id}
\tr \big(P(z,\bar z)\langle z,z\rangle\big) = (n+p+q)P(z,\bar z) 
+\big (\tr P(z,\bar z)\big)\langle z,z\rangle.
\end{equation}
\el

\bpf
By straightforward calculations, we have
\begin{multline}
\tr \big(P(z,\bar z)\langle z,z\rangle\big)
= \sum_{j} \eps_j \frac{\d^2}{\d z_j \d \bar z_j}
\big(P(z,\bar z)\sum_s \eps_s z_s \bar z_s\big)\\
=\big(\sum_{j} \eps_j \frac{\d^2}{\d z_j \d \bar z_j}
P(z,\bar z)\big)\sum_s \eps_s z_s \bar z_s
+\sum_{j} \eps_j P_{z_j}(z,\bar z) \eps_j z_j
+\sum_{j} \eps_j P_{\bar z_j}(z,\bar z) \eps_j \bar z_j
+P(z,\bar z) \sum_j \eps_j^2.
\end{multline}
Using \eqref{hom} we obtain the right-hand side of \eqref{tr-id}
as desired.
\epf

\bpf[Proof of Proposition~\ref{main-decomp} for $s=1$]
We begin by proving the uniqueness of the decomposition \eqref{decomp}.
Let $P(z,\bar z)$ be bihomogeneous in $(z,\bar z)$ of bidegree $(p,q)$ and
suppose that 
\begin{equation}\Label{dec-1}
P(z,\bar z) = Q(z,\bar z) \langle z,z\rangle + R(z,\bar z), \quad \tr R=0.
\end{equation}
Applying $k\ge 1$ times $\tr$ to both sides of \eqref{dec-1}
and using \eqref{tr-id} we obtain, by induction on $k$,
\begin{equation}\Label{k-id}
\tr^k P(z,\bar z) = c_k \tr^{k-1} Q(z,\bar z) + \big( \tr^k Q(z,\bar z)\big)\langle z,z\rangle,
\end{equation}
where $c_k$ are positive integers depending only on $n,p,q$ and
satisfying
\begin{equation}\Label{ck-id}
c_1= n+p+q-2, \quad c_{k+1}=c_k + n+p+q-2k-2.
\end{equation}
Since applying $\tr$ decreases both degrees in $z$ and $\bar z$ by $1$,
one has $\tr^{k_0} Q=0$ for $k_0:=\min (p,q)$.
Hence $\tr^{k_0-1}Q$ is uniquely determined from \eqref{k-id} for $k=k_0$.
Then going backwards through the identities \eqref{k-id}
for $k=k_0-1,k_0-2,\ldots,1$, we see that each $\tr^{k-1}Q$ is uniquely determined
including $\tr^0Q=Q$. Then $R$ is uniquely determined by \eqref{dec-1} proving the uniqueness part of Proposition~\ref{main-decomp} for $s=1$.

To prove the existence, consider the equations
\begin{equation}\Label{k-id1}
\tr^k P(z,\bar z) = c_k Q_{k-1}(z,\bar z) + Q_k(z,\bar z) \langle z,z\rangle, \quad k=1,\ldots,k_0,
\end{equation}
obtained from \eqref{k-id} by replacing each $\tr^k Q(z,\bar z)$
with an indeterminant polynomial $Q_k(z,\bar z)$,
where $Q_{k_0}=0$ for bidegree reason. Hence the last equation for $k=k_0$ reads
$\tr^{k_0}P = c_{k_0} Q_{k_0-1}$,
which we can solve for $Q_{k_0-1}$.
Going backwards through the equations \eqref{k-id1} for $k=k_0-1,k_0-2,\ldots,1$,
as before, we can solve the system \eqref{k-id1} uniquely for $Q_{k_0-2},\ldots,Q_0$.
We claim that 
\begin{equation}\Label{tr-rel}
Q_k=\tr Q_{k-1}, \quad k=1,\ldots,k_0.
\end{equation}
Indeed, \eqref{tr-rel} clearly holds for $k=k_0$ for bidegree reason.
Suppose \eqref{tr-rel} holds for $k>k_1$.
Applying $\tr$ to both sides of \eqref{k-id1} for $k=k_1$, and using \eqref{tr-id} we obtain
\begin{equation}\Label{k-id1'}
\tr^{k_1+1} P(z,\bar z) = c_{k_1} \tr Q_{k_1-1}(z,\bar z) + (n+p+q-2k-2) Q_{k_1}(z,\bar z)  
+ \big(\tr Q_{k_1}(z,\bar z)\big) \langle z,z\rangle, 
\end{equation}
which we compare to \eqref{k-id1} for $k=k_1+1$:
\begin{equation}\Label{k-id1''}
\tr^{k_1+1} P(z,\bar z) = c_{k_1+1} Q_{k_1}(z,\bar z) + \big (\tr Q_{k_1}(z,\bar z)\big) \langle z,z\rangle,
\end{equation}
where we have used \eqref{tr-rel} for $k=k_1+1$.
Using \eqref{ck-id}, we immediately obtain \eqref{tr-rel} for $k=k_1$.
Thus \eqref{tr-rel} holds for all $k$ by induction.
In particular, substituting into \eqref{k-id1} for $k=1$, we obtain
\begin{equation}
\tr P(z,\bar z) = c_1 Q_0(z,\bar z) + \big(\tr Q_0(z,\bar z) \big) \langle z,z\rangle 
= \tr \big( Q_0(z,\bar z) \langle z,z\rangle \big),
\end{equation}
where we have used \eqref{tr-id} and \eqref{ck-id}.
Thus we can take $Q(z,\bar z):=Q_0(z,\bar z)$ and 
$R(z,\bar z):= P(z,\bar z) - Q_0(z,\bar z)\langle z,z\rangle $ to satisfy 
\eqref{dec-1}, proving the existence part of Proposition~\ref{main-decomp} for $s=1$.
\epf

In the Proposition~\ref{main-decomp} in general case we shall use the following lemma.

\bl\Label{trs-lemma}
Assume that a polynomial $P(z,\bar z)$ satisfies $\tr P=0$. Then
\begin{equation}\Label{trs-id}
\tr^s\big(P(z,\bar z)\langle z,z\rangle^{s-1}\big)= 0
\end{equation}
for any $s\ge 1$.
\el

\bpf
 Using \eqref{tr-id}
for $P(z,\bar z)$ replaced with $P(z,\bar z)\langle z,z\rangle^{s-1}$ we obtain
\begin{equation}\Label{trs}
\tr\big(P(z,\bar z)\langle z,z\rangle^{s}\big)
=\tr\big(P(z,\bar z)\langle z,z\rangle^{s-1} \langle z,z\rangle\big)
=c_sP(z,\bar z)\langle z,z\rangle^{s-1} + \big(\tr \big(P(z,\bar z)\langle z,z\rangle^{s-1}\big)\big)
\langle z,z\rangle
\end{equation}
for suitable integer $c_s$ depending on $s$. Replacing $\tr \big(P(z,\bar z)\langle z,z\rangle^{s-1}\big)$
with the right-hand side of \eqref{trs} for $s$ replaced with $s-1$ 
and continuing the process we obtain
\begin{equation}\Label{tr-rel}
\tr\big(P(z,\bar z)\langle z,z\rangle^{s}\big) = c'_sP(z,\bar z)\langle z,z\rangle^{s-1},
\end{equation}
where $c'_s$ is another integer depending on $s$.
Now \eqref{trs-id} can be proved directly by induction on $s$ using \eqref{tr-rel}.
\epf

\bpf[Proof of Proposition~\ref{main-decomp} in general case]
We prove the statement by induction on $s$.
Since it has been proved for $s=1$, it remains to prove the induction step.
Suppose for given $P$ we have a decomposition \eqref{decomp}.
Applying Proposition~\ref{main-decomp} for $s=1$ to $Q$, we also obtain
\begin{equation}
Q(z,\bar z) = Q'(z,\bar z)\langle z,z\rangle  + R'(z,\bar z), \quad \tr R'=0.
\end{equation}
Substituting into \eqref{decomp}, we obtain
\begin{equation}
P(z,\bar z) = Q'(z,\bar z)\langle z,z\rangle^{s+1} +R'(z,\bar z)\langle z,z\rangle^{s} + R(z,\bar z), \quad \tr^sR=0, \quad \tr R'=0.
\end{equation}
Furthermore, $\tr^{s+1}\big(R'(z,\bar z)\langle z,z\rangle^{s}\big)=0$ by Lemma~\ref{trs-lemma}
and therefore we obtain a decomposition
\begin{equation}
P(z,\bar z) = Q'(z,\bar z)\langle z,z\rangle^{s+1} + R''(z,\bar z), \quad \tr^{s+1} R''=0,
\end{equation}
as desired with $R''(z,\bar z):= R'(z,\bar z)\langle z,z\rangle^{s} + R(z,\bar z)$.
This proves the existence part.

Clearly it suffices to prove the uniqueness for $P=0$.
Assume it holds for $s$ and that there is another decomposition
\begin{equation}\Label{2Q}
0 = \2Q(z,\bar z)\langle z,z\rangle^{s+1}  + \2R(z,\bar z), \quad \tr^{s+1}\2R=0.
\end{equation}
By Proposition~\ref{main-decomp} for $s=1$, we can write
\begin{equation}\Label{2R}
\2R = \2Q'(z,\bar z)\langle z,z\rangle + \2R'(z,\bar z), \quad \tr \2R'=0.
\end{equation}
Substitution into \eqref{2Q} yields
\begin{equation}\Label{2Q'}
0 = \big(\2Q(z,\bar z)\langle z,z\rangle^{s}  +  \2Q'(z,\bar z)\big)\langle z,z\rangle + \2R'(z,\bar z), \quad \tr\2R'=0.
\end{equation}
Then the uniqueness for $s=1$ implies $\2R'=0$ and
\begin{equation}\Label{QQ}
\2Q(z,\bar z)\langle z,z\rangle^{s}  +  \2Q'(z,\bar z) =0.
\end{equation}
Applying $\tr^{s+1}$ to both sides of \eqref{2R} we obtain
\begin{equation}
\tr^{s+1} \big( \2Q'(z,\bar z)\langle z,z\rangle\big)=0.
\end{equation}
Now using the identities \eqref{k-id} for $P(z,\bar z):=\2Q'(z,\bar z)\langle z,z\rangle$,
$Q$ replaced by $\2Q'$ and $k=s+1,\ldots, k_0$ with $k_0$ as chosen there, 
we conclude that $\tr^{s}\2Q'=0$.
Then the uniqueness in \eqref{QQ} implies $\2Q=0$.
Hence $\2R=0$ by \eqref{2Q} and
the proof is complete.
\epf

\section{Normalizations}
We now proceed with normalization of the equation for $M'$,
i.e.\ of the terms $\phi'_{klm}(z,\bar z)$.
As noted at the end of \S\ref{weight}, we may assume by induction on the weight
that all ``bad'' terms denoted by dots in \eqref{dots-eq1} -- \eqref{dots-eq3} are fixed.
As explained in  \S\ref{weight}, every coefficient $g_{kl}$ (resp.\ $f_{kl}$)
contributes to ``good" terms corresponding to the lines through $(k,0,l)$, $(0,k,l)$
(resp.\ $(k,1,l)$, $(1,k,l)$) in the same direction $(1,1,-1)$.

\subsection{Normalization for $k\ge 2$}\Label{g2}
We begin by analysing the line through $(k,0,l)$ for $k\ge 2$,
corresponding to the multi-degrees $(k+m,m,l-m)$, $0\le m\le l$,
for which we have the equation \eqref{dots-eq1}.
We rewrite this equation as
\begin{multline}\Label{dots-eq1'}
\phi'_{k+m,m,l-m}(z,\bar z) =\frac{1}{2i}{l\choose m} g_{kl}(z)(i\langle z,z\rangle)^m
- {l-1\choose m-1}\langle f_{k+1,l-1}(z), z\rangle (i\langle z,z\rangle)^{m-1}
+  \ldots,
\end{multline}
Where we have included the given term $\phi_{k+m,m,l-m}(z,\bar z)$ in the dots.
Our goal is to write normalization conditions for $\phi'_{k+m,m,l-m}(z,\bar z)$
that uniquely determine $g_{kl}$ and $f_{k+1,l-1}$.
If $m=0$, the term with $f_{k+1,l-1}$ is not present.
Hence we have to consider an identity \eqref{dots-eq1'} with $m\ge 1$.
Then the sum of the terms involving $g_{kl}$ and $f_{k+1,l-1}$ is a multiple 
of $\langle z,z\rangle^{m-1}$. Thus, by varying $g_{kl}$ and $f_{k+1,l-1}$,
we may expect to normalize $\phi'_{k+m,m,l-2m}(z,\bar z)$ to be in the complement
of the space of polynomials of the form $P(z,\bar z)\langle z,z\rangle^{m-1}$.
The suitable decomposition is given by Proposition~\ref{main-decomp}:
\begin{equation}\Label{phi'-dec}
\phi'_{k+m,m,l-m}(z,\bar z) = Q(z,\bar z)\langle z,z\rangle^{m-1} + R(z,\bar z), \quad \tr^{m-1}R=0.
\end{equation}
Using similar decompositions for other terms in \eqref{dots-eq1'} and 
equating the factors of $\langle z,z\rangle^{m-1}$, we obtain
\begin{equation}\Label{m}
Q(z,\bar z) =\frac{1}{2i}{l\choose m} g_{kl}(z)i^m\langle z,z\rangle
- {l-1\choose m-1}\langle f_{k+1,l-1}(z), z\rangle i^{m-1}
+  \ldots,
\end{equation}
where, as before, the dots stand for the terms that have been fixed.
Since $f_{k+1,l-1}(z)$ is free and the form $\langle z,z\rangle$ is nondegenerate,
$\langle f_{k+1,l-1}(z), z\rangle$ is a free bihomogeneous polynomial of the corresponding bidegree. Hence we can normalize $Q$ to be zero. In view of \eqref{phi'-dec}, $Q=0$
is equivalent to the normalization condition
\begin{equation}\Label{norm1}
\tr^{m-1}\phi'_{k+m,m,l-m}=0.
\end{equation}
Putting $Q=0$ in \eqref{m},
we can now uniquely solve this equation for $\langle f_{k+1,l-1}(z), z\rangle$ in the form
\begin{equation}\Label{f-solv}
\langle f_{k+1,l-1}(z), z\rangle = \frac{1}{2}{l\choose m}{l-1\choose m-1}^{-1} g_{kl}(z) \langle z,z\rangle+ \ldots,
\end{equation}
from where $f_{k+1,l-1}(z)$ is uniquely determined.
It remains to determine $g_{kl}$, for which we substitute \eqref{f-solv}
into an identity \eqref{dots-eq1'} with $m$ replaced by $m'\ne m$:
\begin{multline}\Label{dots-eq1''}
\phi'_{k+m',m',l-m'}(z,\bar z) =\frac{1}{2i}\Big({l\choose m'} 
- {l-1\choose m'-1}{l\choose m}{l-1\choose m-1}^{-1}\Big) g_{kl}(z)(i\langle z,z\rangle)^{m'} +  \ldots,
\end{multline}
where we have assumed $m'\ge1$. If $m'=0$, the term with $f_{k+1,l-1}$ does not occur
and hence no substitution is needed. In the latter case, the coefficient in front of $g_{kl}(z)$
is clearly nonzero. Otherwise, for $m'\ge1$, that coefficient is equal, up to a nonzero factor, to the determinant
\begin{equation}\Label{det}
\left|
\begin{matrix}
{l\choose m}&-{l-1\choose m-1}\\
{l\choose m'}&-{l-1\choose m'-1}
\end{matrix}
\right|
=
\left|
\begin{matrix}
{l!\over m!(l-m)!}&{(l-1)!\over (m-1)!(l-m)!}\\
{l!\over m'!(l-m')!}&{(l-1)!\over (m'-1)!(l-m')!}
\end{matrix}
\right|
=
\frac{l!(l-1)!}{m!(l-m)!m'!(l-m')!}
\left|
\begin{matrix}
{1}&{m}\\
{1}&{m'}
\end{matrix}
\right|\ne0.
\end{equation}
Hence the coefficient in front of $g_{kl}(z)$ in \eqref{dots-eq1''}
is nonzero in any case.
Therefore, using Proposition~\ref{main-decomp} as above
we see that we can obtain the normalization condition
\begin{equation}\Label{norm1'}
\tr^{m'}\phi'_{k+m',m',l-m'}=0,
\end{equation}
which determines uniquely $g_{kl}(z)$
and hence $f_{k+1,l-1}(z)$ in view of \eqref{f-solv}.

Summarizing, for each $m\ge1$ and $m'\ne m$, we obtain
the normalization conditions
\begin{equation}\Label{norm1''}
\tr^{m-1}\phi'_{k+m,m,l-m}=0,\quad \tr^{m'}\phi'_{k+m',m',l-m'}=0,
\end{equation}
which determine uniquely $g_{kl}(z)$ and $f_{k+1,l-1}(z)$.
Such a choice of $m$ and $m'$ is always possible unless $l=0$.
In the latter case, the coefficient $f_{k+1,-1}(z)=0$ is not present
and $g_{k0}(z)$ is uniquely determined by
the normalization condition $\phi'_{k00}=0$ corresponding to $m'=0$.
That is, for $l=0$, we only have the second condition in \eqref{norm1''} with $m'=0$.

Since the degree of $\phi'_{k+m,m,l-m}$ is $k+l+m$,
we obtain the lowest possible degrees in \eqref{norm1''} for $m'=0$, $m=1$,
which corresponds to the normalization
\begin{equation}
\phi'_{k+1,1,l-1}=0, \quad \phi'_{k0l}=0,
\end{equation}
which is the part of the Chern-Moser normal form \cite{CM}.

\subsection{Normalization for $k=1$}
We next analyze the line through $(1,0,l)$ corresponding to the multi-degrees $(m+1,m,l-m)$,
$0\le m\le l$, for which we have the equation \eqref{dots-eq2}. As before, we rewrite this equation as
\begin{multline}\Label{dots-eq2'}
\phi'_{m+1,m,l-m}(z,\bar z) = \frac{1}{2i}{l\choose m} g_{1l}(z)(i\langle z,z\rangle)^m  \\
- {l-1\choose m-1}\langle f_{2,l-1}(z), z\rangle (i\langle z,z\rangle)^{m-1}
-{l-1\choose m-1}\langle z, f_{0l}\rangle  (-i\langle z,z\rangle)^m + \ldots.
\end{multline}
Arguing as before, we consider $m\ge 1$ and decompose
\begin{equation}
\phi'_{m+1,m,l-m}(z,\bar z) = Q(z,\bar z)\langle z,z\rangle^{m-1} + R(z,\bar z), \quad
\tr^{m-1} R=0.
\end{equation}
Decomposing similarly the other terms in \eqref{dots-eq2'} and equating the factors of
$\langle z,z\rangle^{m-1}$, we obtain
\begin{multline}\Label{dots-eq2'}
Q(z,\bar z) = \frac{1}{2i}{l\choose m} g_{1l}(z)i^m\langle z,z\rangle  \\
- {l-1\choose m-1}\langle f_{2,l-1}(z), z\rangle  i^{m-1}
-{l-1\choose m-1}\langle z, f_{0l}\rangle  (-i)^m\langle z,z\rangle + \ldots.
\end{multline}
Since $f_{2,l-1}(z)$ is free and the form $\langle z,z\rangle$ is nondegenerate,
$\langle f_{2,l-1}(z), z\rangle$ is also free and hence we can choose it suitably to
obtain $Q=0$, which is equivalent to the normalization condition
\begin{equation}
\tr^{m-1}\phi'_{m+1,m,l-m}=0.
\end{equation}
Putting $Q=0$ into \eqref{dots-eq2'}, we solve it uniquely for
$\langle f_{2,l-1}(z), z\rangle$, which, in turn, 
determines uniquely $f_{2,l-1}(z)$.
Arguing as in \S\ref{g2}, we substitute the obtained  expression for 
$\langle f_{2,l-1}(z), z\rangle$ into identities \eqref{dots-eq2'}
with $m$ replaced by $m'\ne m$. The result can be written as
\begin{equation}
\phi'_{m'+1,m',l-m'}(z,\bar z) = \big(c_{m'}g_{1l}(z) + d_{m'}\langle z, f_{0l}\rangle\big)
\langle z,z\rangle^{m'} + \ldots
\end{equation}
with suitable coefficients $c_{m'}$, $d_{m'}$.
Then as in \S\ref{g2} we see that we can obtain the normalization
\begin{equation}\Label{norm2}
\tr^{m'}\phi'_{m'+1,m',l-m'}=0, \quad \tr^{m''}\phi'_{m''+1,m'',l-m''}=0
\end{equation}
for any pair $(m',m'')$ such that
\begin{equation}\Label{nonzero}
\left|
\begin{matrix}
c_{m'} & d_{m'}\\
c_{m''} & d_{m''}
\end{matrix}
\right|\ne 0.
\end{equation}
Furthermore, assuming \eqref{nonzero}, both $g_{1l}(z)$ and $f_{0l}$
are uniquely determined by \eqref{norm2}, which, in turn, determine 
$f_{2,l-1}(z)$ in view of \eqref{dots-eq2'} with $Q=0$.
Thus it remains to choose $m',m''$ satisfying \eqref{nonzero}.
Inspecting the construction of the coefficients $c_{m'}$, $d_{m'}$,
it is straightforward to see that the determinant in \eqref{nonzero} is equal,
up to a nonzero multiple, to the determinant
\begin{equation}\Label{nonzero''}
\left|
\begin{matrix}
\frac{1}{2i}{l\choose m}i^m & -{l-1\choose m-1}i^{m-1} & - {l-1\choose m-1}(-i)^m\\
\frac{1}{2i}{l\choose m'}i^{m'} & -{l-1\choose m'-1}i^{m'-1} & - {l-1\choose m'-1}(-i)^{m'}\\
\frac{1}{2i}{l\choose m''}i^{m''} & -{l-1\choose m''-1}i^{m''-1} & - {l-1\choose m''-1}(-i)^{m''}\\
\end{matrix}
\right|
\end{equation}
consisting of the coefficients in \eqref{dots-eq2'}.
In fact one can see that the matrix in \eqref{nonzero}
is obtained as a block in \eqref{nonzero'} by elementary row operations.
The determinant \eqref{nonzero''} is equal, up to a nonzero multiple, to
 \begin{equation}\Label{nonzero'}
\left|
\begin{matrix}
1 & m & (-1)^m m\\
1 & m' & (-1)^{m'} m'\\
1 & m'' & (-1)^{m''} m''\\
\end{matrix}
\right|=
\left|
\begin{matrix}
m'-m & (-1)^{m'} m'-(-1)^m m\\
m''-m & (-1)^{m''} m''-(-1)^m m\\
\end{matrix}
\right|=
\pm\left|
\begin{matrix}
m'-m & (-1)^{m'-m} m'- m\\
m''-m & (-1)^{m''-m} m''- m\\
\end{matrix}
\right|.
\end{equation}
In case both $m'-m$ and $m''-m$ are odd,
the latter determinant is equal to
\begin{equation}
(m'-m)(-m''-m) - (m''-m)(-m'-m) = 2m''m-2m'm=2m(m''-m),
\end{equation}
which is nonzero by the construction.
In case $m'-m$ is even and $m''-m$ is odd,
the last determinant in \eqref{nonzero'} is equal, up to a sign, to
\begin{equation}
(m'-m)(-m''-m) - (m''-m)(m'-m) = 2m''m-2m''m' = 2m''(m-m'),
\end{equation}
which is nonzero provided $m''\ne0$.
Similarly, in case $m'-m$ is odd and $m''-m$ is even,
the determinant is nonzero provided $m'\ne0$.
On the other hand, if $m''=0$, the determinant \eqref{nonzero'} is nonzero if $m'-m$ is odd
and similarly, if $m'=0$, the determinant \eqref{nonzero'} is nonzero if $m''-m$ is odd.
In all other cases, the determinant is zero.

Summarizing we conclude that \eqref{nonzero'} is nonzero and hence \eqref{nonzero} holds
whenever $m,m',m''$ are disjoint and the nonzero ones among them are not of the same parity.
For such a choice of $m,m',m''$ with $m\ge1$, we have the normalization conditions
\begin{equation}\Label{norm2}
\tr^{m-1}\phi'_{m+1,m,l-m}=0, \quad \tr^{m'}\phi'_{m'+1,m',l-m'}=0, \quad \tr^{m''}\phi'_{m''+1,m'',l-m''}=0,
\end{equation}
that determine uniquely $g_{1l}(z)$, $f_{2,l-1}(z)$ and $f_{0l}$.
Such a choice of $m,m',m''$ is always possible unless $l\in\{0,1\}$.
In case $l=0$, all terms in \eqref{dots-eq2'} are already zero.
If $l=1$, the only choice is $m=1$ and $m'=0$ leaving no space for $m''$.
In that case we still obtain the normalization
\begin{equation}\Label{norm2''}
\phi'_{101}=0, \quad \phi'_{210}=0,
\end{equation}
which, however, does not determine $g_{11}(z)$, $f_{20}(z)$ and  $f_{01}$
uniquely. Instead, we regard $f_{01}$ as a free parameter
and then \eqref{norm2''} determine uniquely $g_{11}(z)$ and $f_{20}(z)$.
The free parameter $f_{01}$ corresponds to the choice of $a\in\C^n$
in the following group of automorphisms of the quadric $\Im w=\langle z,z\rangle$:
\begin{equation}
(z,w)\mapsto \frac{(z+aw,w)}{1-2i\langle z,a\rangle-i\langle a,a\rangle w}.
\end{equation}

To obtain the lowest possible degrees in \eqref{norm2},
we have to choose $m=1$, $m'=0$, $m''=0$, in which case \eqref{norm2} gives
\begin{equation}
\phi'_{10l}=0, \quad \phi'_{2,1,l-1}=0, \quad \tr^2\phi'_{3,2,l-2}=0,
\end{equation}
which is precisely a part of the Chern-Moser normal form \cite{CM}.

\subsection{Normalization for $k=0$}
It remains to analyze the line through $(0,0,l)$ corresponding
to the multi-degrees $(m,m,l-m)$, $0\le m\le l$, for which we have the equation
\eqref{dots-eq3}. Similarly to the above, we rewrite this equation as
\begin{multline}\Label{dots-eq3'}
\phi'_{m,m,l-m}(z,\bar z)=
{l\choose m}\Im \big(g_{0l}(i\langle z,z\rangle)^m\big) -
2{l-1\choose m-1}\Re\big(\langle f_{1,l-1}(z), z\rangle (i\langle z,z\rangle)^{m-1}\big) + \ldots.
\end{multline}
In view of the presence of the power $i^m$ inside real and imaginary parts,
we have to deal separately with $m$ being even and odd.

We first assume $m$ to be even. In that case \eqref{dots-eq3'} can be rewritten as
\begin{multline}\Label{dots-eq3''}
\phi'_{m,m,l-m}(z,\bar z)=
{l\choose m}(\Im g_{0l})(i\langle z,z\rangle)^m -
2i{l-1\choose m-1}(\Im\langle f_{1,l-1}(z), z\rangle) (i\langle z,z\rangle)^{m-1}+ \ldots.
\end{multline}
Argueing as above we can obtain the normalization
\begin{equation}
\tr^{m-1}\phi'_{m,m,l-m}=0,
\end{equation}
implying an equation for $f_{1,l-1}(z)$ and $g_{0l}$
which can be solved for $\Im\langle f_{1,l-1}(z), z\rangle$ in the form
\begin{equation}
\Im\langle f_{1,l-1}(z), z\rangle = \frac12 {l-1\choose m-1}^{-1} {l\choose m} (\Im g_{0l})
\langle z,z \rangle + \ldots.
\end{equation}
The latter expression is to be substituted into another equation \eqref{dots-eq3''}
with $m$ replaced by $m'$ (still even). 
As above, we obtain the normalization
\begin{equation}\Label{norm3}
\tr^{m-1}\phi'_{m,m,l-m}=0, \quad  \tr^{m'}\phi'_{m',m',l-m'}=0,
\end{equation}
provided
\begin{equation}
\left|
\begin{matrix}
{l\choose m}&-2{l-1\choose m-1}\\
{l\choose m'}&-2{l-1\choose m'-1}
\end{matrix}
\right|\ne 0,
\end{equation}
which always holds in view of \eqref{det}.
Hence $\Im g_{0l}$ is uniquely determined by \eqref{norm3}
and therefore also $\Im\langle f_{1,l-1}(z), z\rangle$.

Now assume $m$ is odd. In that case \eqref{dots-eq3'} becomes
\begin{multline}\Label{dots-eq3''}
\phi'_{m,m,l-m}(z,\bar z)=
{l\choose m}(\Re g_{0l})i^{m-1}\langle z,z\rangle^m -
2{l-1\choose m-1}(\Re\langle f_{1,l-1}(z), z\rangle) (i\langle z,z\rangle)^{m-1}+ \ldots.
\end{multline}
Then the above argument yields the normalization
\begin{equation}\Label{norm3'}
\tr^{m-1}\phi'_{m,m,l-m}=0, \quad  \tr^{m'}\phi'_{m',m',l-m'}=0,
\end{equation}
where this time both $m$ and $m'$ are odd.
This time \eqref{norm3'} determines uniquely both 
$\Re g_{0l}$ and $\Re\langle f_{1,l-1}(z), z\rangle$.

Summarizing, we obtain the normalization 
\begin{equation}\Label{norm3'}
\tr^{m-1}\phi'_{m,m,l-m}=0, \quad  \tr^{m'}\phi'_{m',m',l-m'}=0,\quad
\tr^{\2m-1}\phi'_{\2m,\2m,l-\2m}=0, \quad  \tr^{\2m'}\phi'_{\2m',\2m',l-\2m'}=0,
\end{equation}
where $m\ge1$, $m'\ne m$, $\2m'\ne \2m$, both $m,m'$ are even and both $\2m,\2m'$ are odd
(note that automatically $\2m\ge1$).
The choice of such numbers $m,m',\2m,\2m'$ is always possible unless $l\in\{0,1,2\}$.
If $l=0$, all terms in \eqref{dots-eq3'} are zero.
If $l=1$, we obtain a valid identity since we have assumed $g_{01}=1$, $f_{01}=\id$.
Finally, for $l=2$, we can choose $m=2$, $m'=0$ for the even part
and $\2m=1$ for the odd part but there is no place for $\2m'$.
We obtain the normalization
\begin{equation}
\phi'_{002}=0, \quad \phi'_{111}=0, \quad \tr\phi'_{220}=0,
\end{equation}
which determines uniquely $\Im g_{02}(z)$ and $f_{11}(z)$,
whereas $\Re g_{02}(z)$ is a free parameter.
The latter corresponds to the choice of $r\in\R$
in the following group of automorphisms of the quadric 
$\Im w=\langle z,z\rangle$:
\begin{equation}
(z,w)\mapsto \frac{(z,w)}{1-rw}.
\end{equation}

To obtain the lowest possible degrees in \eqref{norm3'}, 
we can choose $m=2$, $m'=0$, $\2m=1$, $\2m'=3$,
which leads to the normalization
\begin{equation}\Label{cm}
\phi'_{00l}=0, \quad \phi'_{1,1,l-1}=0, \quad \tr\phi'_{2,2,l-2}=0, \quad \tr^3\phi'_{3,3,l-3}=0,
\end{equation}
which is precisely a part of Chern-Moser normal form \cite{CM}.
However, there is another choice in the lowest degree, namely
$m=2$, $m'=0$, $\2m=3$, $\2m'=1$, in which case
the normalization reads
\begin{equation}\Label{cm'}
\phi'_{00l}=0, \quad \tr\phi'_{1,1,l-1}=0, \quad \tr\phi'_{2,2,l-2}=0, \quad \tr^2\phi'_{3,3,l-3}=0.
\end{equation}
Comparing \eqref{cm} and \eqref{cm'}
we can say that the Chern-Moser normalization \eqref{cm}
has more equations for $\phi'_{1,1,l-1}$ and less equations for $\phi'_{3,3,l-3}$.
In a sense, the Chern-Moser normalization corresponds to the maximum conditions
in the lowest possible degree.

Summarizing the results of this section, we obtain the proof of Theorem~\ref{main}.

\end{document}